\documentclass[twocolumn,amsthm]{autart}

\usepackage[utf8]{inputenc}
\usepackage[T1]{fontenc}
\usepackage{lmodern}

\usepackage{mathtools}
  \mathtoolsset{showonlyrefs, centercolon, mathic}
\usepackage{amssymb}
\usepackage{dsfont}
\usepackage{nicefrac}
\usepackage[np,autolanguage]{numprint}
\usepackage{upgreek}

\usepackage[usenames, dvipsnames, svgnames]{xcolor}
\usepackage{tikz}
  \usetikzlibrary{arrows}
    \tikzset{>=stealth}
  \usetikzlibrary{external}
    \tikzset{external/only named}
    \tikzset{external/mode=list and make}
  \usetikzlibrary{plotmarks}
  \usetikzlibrary{positioning}

\usepackage{pgfplots}
  \usepgfplotslibrary{groupplots}
  \usepgfplotslibrary{statistics}
  \pgfplotsset{compat=newest}
  \pgfplotsset{
    log ticks with fixed point,
    clip marker paths=true,
  }

\usepackage{enumerate}
\usepackage{dblfloatfix}
\usepackage{paralist}
\usepackage[multidot]{grffile}
\usepackage{needspace}
\usepackage{xifthen}
\usepackage{todonotes}

\usepackage{natbib}

\usepackage[colorlinks=true, linkcolor=black, citecolor=black, 
            urlcolor=blue]{hyperref}

\newcommand{\inv}{^{-1}} 

\newcommand{\downto}{\downarrow} 

\newcommand{\dd}{\ensuremath{\mathrm{d}}} 
\newcommand{\reals}{\ensuremath{\mathds{R}}} 
\newcommand{\realspos}{\ensuremath{\reals_{>0}}} 
\DeclareMathOperator{\diverg}{div} 
\DeclareMathOperator{\support}{supp} 

\newcommand{\borelsa}[1]{\mathcal{B}_{#1}} 

\newcommand{\paren}[1]{\ensuremath{\left(#1\right)}}

\newcommand{\fparen}[1]{\ensuremath{\!\left(#1\right)}}

\newcommand{\norm}[2][]{\ensuremath{\left\lVert#2\right\rVert}_{#1}}

\newcommand{\normnm}[2][]{\ensuremath{\lVert#2\rVert}_{#1}}

\newcommand{\Prob}[2][]{%
    \ensuremath{P\ifthenelse{\isempty{#1}}{}{_{#1}\!}\fparen{#2}}%
}

\newcommand{\pr}[2][]{%
  \ensuremath{p\ifthenelse{\isempty{#1}}{\!}{_{#1}\negmedspace}\paren{#2}}%
}

\newcommand{\snorm}[1]{%
  \left\lvert\kern-0.3ex\left\lvert\kern-0.3ex\left\lvert
  #1
  \right\rvert\kern-0.3ex\right\rvert\kern-0.3ex\right\rvert
}

\newcommand{\snormnm}[1]{%
  \lvert\kern-0.3ex\lvert\kern-0.3ex\lvert
  #1
  \rvert\kern-0.3ex\rvert\kern-0.3ex\rvert
}

\newcommand{\E}[2][]{\ensuremath{\mathrm{E}_{#1}\negmedspace\left[#2\right]}}

\newcommand{\normald}[1]{\ensuremath{\mathcal{N}\!\paren{#1}}}

\newcommand{\cond}{\ensuremath{\,\middle\vert\,}}

\newcommand{\inputtikzpicture}[1]{%
  \includegraphics{#1.tikz.pdf}%
}

\newcommand{\nx}{m}
\newcommand{\nz}{n}
\newcommand{\npar}{q}

\newcommand{\cmspace}{{\mathcal H}}
\newcommand{\tspace}{{\mathcal T}}
\newcommand{\yspace}{{\mathcal Y}}
\newcommand{\filti}{\mathcal F_0}

\newcommand{\sigmad}{\sigma_{\mathrm d}}
\newcommand{\sigmai}{\sigma_{\mathrm 0}}
\newcommand{\sigmao}{\sigma_{\mathrm o}}
\newcommand{\sigmap}{\sigma_{\mathrm \theta}}
\newcommand{\sigmar}{\sigma_{\mathrm r}}
\newcommand{\sigmay}{\sigma_{\mathrm y}}
\newcommand{\Sigmay}{\Sigma_{\mathrm y}}
\newcommand{\ts}{t_{\mathrm s}}
\newcommand{\po}{p_{\mathrm o}}
\newcommand{\lenerg}{\ell_{\mathrm e}}

\begin{document}

\begin{frontmatter}
  \title{%
    Joint Maximum \emph{a Posteriori} \ State Path and Parameter Estimation
    in Stochastic Differential Equations\thanksref{footnoteinfo}
  }

  \thanks[footnoteinfo]{%
    This work has been financially supported by the Brazilian agencies 
    CAPES, CNPq and FAPEMIG.
  }
  
  \author[PPGMEC]{Dimas Abreu Archanjo Dutra}\ead{dimasad@ufmg.br},
  \author[PPGEE]{Bruno Otávio Soares Teixeira}\ead{brunoot@ufmg.br},
  \author[PPGEE]{Luis Antonio Aguirre}\ead{aguirre@cpdee.ufmg.br}.

  \address[PPGMEC]{%
    Programa de Pós-Graduação em Engenharia Mecânica ---
    Universidade Federal de Minas Gerais ---
    Belo Horizonte, MG, Brasil
  }
  \address[PPGEE]{%
    Programa de Pós-Graduação em Engenharia Elétrica ---
    Universidade Federal de Minas Gerais ---
    Belo Horizonte, MG, Brasil
  }
  
  
  \begin{abstract}
    In this article, we introduce the joint maximum \emph{a posteriori}
    state path and parameter estimator (JME) for continuous-time systems
    described by stochastic differential equations (SDEs).
    This estimator can be applied to nonlinear systems with
    discrete-time (sampled) measurements with a wide range of measurement
    distributions.
    We also show that the minimum-energy state path and parameter estimator
    (MEE) obtains the joint maximum \emph{a posteriori} noise path,
    initial conditions, and parameters.
    These estimators are demonstrated in simulated experiments, in which they
    are compared to the prediction error method (PEM) using the unscented Kalman
    filter and smoother.
    The experiments show that the MEE is biased for the damping parameters of
    the drift function.
    Furthermore, for robust estimation in the presence of outliers, the
    JME attains lower state estimation errors than the PEM.
  \end{abstract}
\end{frontmatter}

\section{Introduction}
\label{sec:introduction}

In the context of discrete-time systems, maximum \emph{a posteriori} (MAP) 
state path estimators have recently emerged as a powerful alternative to Kalman
filters and smoothers due to their robustness properties and applicability to
a larger class of models 
\citep{
  bell2009icn, aravkin2011llr, aravkin2012rtf, aravkin2012sct, aravkin2013spe,
  farahmand2011drs, dutra2012jmaps, monin2013mte%
}.
A wide variety of phenomena of engineering interest are continuous-time in
nature and can be modeled by stochastic differential equations (SDEs).
For this class of models, the MAP state-path estimator is built upon the
Onsager--Machlup functional and is the solution to an optimal control problem
\citep{zeitouni1987map, aihara1999mln, aihara1999mem}.

To evaluate if a discretization of a variational optimization problem
is \emph{consistent}, the concept of hypographical convergence is used
\citep[cf.][]{polak2011rof}.
If a sequence of discretized problems hypo-converges to a variational problem,
then the discretized optima converge to the variational optima.
In a previous work \citep{dutra2014map}, we showed that the discrete-time MAP
state path estimator applied to trapezoidally discretized continuous-time 
systems converges hypographically to the MAP state path estimator of the
continuous systems, as the discretization step vanishes.
However, when the Euler discretization is used instead---the most widespread 
approach---the discretized estimator hypo-converges to the minimum-energy
estimator, whose estimates were proved to be MAP \emph{noise} paths.
This implies that the discretized MAP estimates have a different interpretation
depending on the discretization scheme used.

In this work, we present the extension of the estimators of \citet{dutra2014map}
for joint state path and parameter estimation.
We introduce the joint MAP state path and parameter estimator
(JME) for continuous-time systems and
also show that the joint minimum-energy state path and parameter estimator
(MEE) corresponds to the joint MAP \emph{noise} path, initial state and
parameter estimator.
The JME and MEE are also the hypographical limits of the 
trapezoidally- and Euler-discretized joint state path and parameter estimators
\citep[Chap.~3]{dutra2014thesis}, respectively.

The fact that the parameter is estimated as a single vector instead of a
time-varying augmented state places the JME and MEE in a similar niche to
the Kalman-filter--based prediction error method (PEM)
\citep{kristensen2004pes}, to which it is compared.
Furthermore, if the joint state path and parameter posterior distribution is
unimodal and approximately symmetric, the JME estimates should be close to 
the marginal MAP parameter estimates.
Similarly, even when the parameters are not of interest, the JME can be used
as a state-path estimator under parametric uncertainty.

The merit function of these estimators admits a tractable expression for a
wide range of nonlinear systems, lending them a wider applicability
than Kalman-filter--based estimators.
In particular, it is possible to use heavy-tailed measurement distributions
which confer robustness against outliers
\citep{
  aravkin2011llr, aravkin2012rtf, aravkin2012sct, aravkin2013spe,
  aravkin2012nrs, farahmand2011drs, dutra2012jmaps%
}. 
The resulting estimators can be seen as an extension of Huber's M-estimators
\citep[Sec.~3.2]{huber2009rs} to the smoothing problem.
M-estimates of location parameters using heavy-tailed distributions
can be interpreted as implicit weighted means, with low weights assigned to
outlying observations.
This approach ``combines conceptual simplicity with generality, 
since it can be applied to a wide range of settings'' 
\citep[p.~882]{lange1989rsm}.
A competing approach to robust estimation is to consider a family of
distributions in the neighbourhood of a nominal guess and design filters or 
smoothers
which guarantee the best behaviour in the worst-case scenario, i.e., minimax
estimators \citep{levy2013rss, zorzi2016rkf}.

The remainder of this article is organized as follows.
In Sec.~\ref{sec:problem_def} we define the problem being tackled and common
notation and variables.
In Sec.~\ref{sec:estimators} the fictitious densities and the MAP estimators
are defined and presented.
The simulated example applications are presented in Sec.~\ref{sec:sim} and
conclusions and future work in Sec.~\ref{sec:conclusions}.

\section{Problem definition}
\label{sec:problem_def}

In what follows, $(\Omega,\mathcal{F}, P)$ is a standard probability space on
which all random variables are defined.
Random variables will be denoted by uppercase letters and their values by 
lowercase, so that if $Y\colon\Omega\to\mathcal{Y}$ is a $\mathcal Y$-valued
random variable, $y\in\mathcal Y$ will denote specific values it might take.
The same applies to stochastic processes.
The dependency on the random outcome $\omega\in\Omega$ will be omitted when
unambiguous, to simplify the notation.
For a random variable $\Theta$, $\support(P_\Theta)$ denotes the topological
support of its induced measure.
The time argument of functions may also be written as subscripts for
compactness.

Let $X$ and $Z$ be $\reals^\nx$- and $\reals^\nz$-valued stochastic processes,
respectively, representing the state of a system over the experiment interval
$\tspace:=[0, T]$ and satisfying the following system of SDEs:
\mathtoolsset{showonlyrefs=false}%
\begin{subequations}
  \label{eq:xz_sde}
  \begin{align}
    \label{eq:x_sde}
    \dd X_t &= f(t, X_t, Z_t, \Theta)\,\dd t + G \,\dd W_t,\\
    \label{eq:z_sde}
    \dd Z_t &= h(t, X_t, Z_t, \Theta)\,\dd t,
  \end{align}
\end{subequations}
\mathtoolsset{showonlyrefs=true}%
where $f$ and $h$ are the drift functions, the $\reals^\npar$-valued random
variable $\Theta$ is the unknown parameter vector, the full rank
$G\in\reals^{\nx\times\nx}$ is the diffusion matrix, and $W$ is an
$\nx$-dimensional Wiener process.
This division of the state in two parts, $X$ and $Z$, is done to cover systems
in which the evolution of some state variables is not under directly influence
of noise.

Consider, in addition, that some $\yspace$-valued random variable $Y$ is
observed.
We assume that the conditional distribution of $Y$, given
$X=x$, $Z=z$ and $\Theta=\theta$, is absolutely continuous and admits a
density $\psi$ with respect to a measure $\nu$ on the measurable space
$(\yspace, \borelsa{\yspace})$,
i.e., for all $\mathbb B\in\borelsa{\yspace}$,
\begin{equation}
  P_Y (\mathbb B \,|\, X=x, Z=z, \Theta=\theta) = 
  \int_{\mathbb B} \psi\fparen{y\cond x,z,\theta}\dd \nu(y).
\end{equation}

In this paper, we derive the joint MAP estimator for $X$, $Z_0$ and $\Theta$,
given $y\in\yspace$.
Note that, conditioned on that, the whole $Z$ path is also uniquely defined.
We also prove that the minimum-energy estimator is the joint MAP estimator for
$W$, $X_0$, $Z_0$ and $\Theta$.
We use the abstract definitions of mode and the MAP estimator of
\citet[Defns.~1 and~2]{dutra2014map}, which cover random variables over
infinite-dimensional spaces such as state paths of continuous-time systems.
These definitions can be better understood using the concept of a
\emph{fictitious density}, which we introduce formally in the definition
below.
Similar terminology was applied to the Onsager--Machlup functional in
this context, it was described as an
\emph{ideal} density with respect to a \emph{fictitious} uniform measure 
by \citet[p.~433]{takahashi1981pf} and as a \emph{fictitious density} by
\citet[p.~1037]{zeitouni1989omf}.

\begin{defn}[{\citealt[Defn.~2.4]{dutra2014thesis}}]
  \label{th:fict_dens}
  Let $A$ be an $\mathcal A$-valued random variable, where $(\mathcal A, d)$
  is a metric space.
  The function $\zeta\colon\mathcal A\to\reals$ is a \emph{fictitious density}
  if $\zeta(a)>0$ for at least one $a\in\mathcal A$ and
  there exists $\xi\colon\realspos\to\realspos$ such that
  \begin{equation}
    \label{eq:fict_dens_defn}
    \lim_{\epsilon\downto 0} 
    \frac{P\fparen{\big.d(A, a)< \epsilon}}{\xi(\epsilon)} = \zeta(a)
    \qquad\text{for all }a\in\mathcal A.
  \end{equation}
\end{defn}

The fictitious density can be understood as a density with respect to a metric.
It quantifies the concentration of probability in the neighbourhood of a point.
When, for some $a',a''\in\mathcal A$, the fictitious density 
$\zeta(a')>\zeta(a'')$, then the $\epsilon$-balls around $a'$ have a larger
probability than those around $a''$, for all sufficiently small $\epsilon$.
This means that the MAP estimates according to \citet[Defn.~2]{dutra2014map}
are the maxima of the posterior fictitious density.
It should be noted that for Euclidean spaces any fictitious density
is proportional to the probability density function in the usual sense.
We now show the application of these concepts to the state paths and
parameters of the system described by \eqref{eq:xz_sde}.

\section{MAP and minimum energy estimators}
\label{sec:estimators}

The following assumptions will be made on the system's probabilistic and 
dynamical model.

\begin{assum}
  \label{th:assum}
  \hfill
  \begin{enumerate}[a.]
  \item \label{it:prior_dens_meas}%
    The initial states $X_0$, $Z_0$ and the parameter vector $\Theta$ are
    $\filti$-measurable and admit a continuous joint prior density \(\pi\) with
    respect to the Lebesgue measure.
  \item \label{it:assum_uniform_cont}
    The functions $f$ and $h$ are uniformly continuous with respect to
    all their arguments for $\theta\in\support(P_\Theta)$.
  \item \label{it:lipschitz}
    For all fixed $\theta\in\support(P_\Theta)$, the functions $f$ and $h$
    are Lipschitz continuous with respect to their second and third arguments
    $x$ and $z$, uniformly over their first argument $t$.
  \item \label{it:fc2}
    For all fixed $\theta\in\support(P_\Theta)$, the function $f$ is twice 
    differentiable with respect to its second argument $x$ and differentiable
    with respect to its first and third arguments $t$ and $z$.
    Furthermore, its first and second derivatives mentioned above are continuous
    with respect to all arguments, for all $\theta\in\support(P_\Theta)$.
  \item \label{it:novikov}
    The system is such that
    \begin{equation}
      \label{eq:novikov}
      \E{
        \exp\fparen{
          \int_0^T \norm{G\inv f(t, X_t, Z_t, \Theta)}^2 \dd t
        }
      } < \infty.
    \end{equation}
  \item The measurement likelihood $\psi$ is continuous with respect to
    the given $x$, $z$ and $\theta$.
  \item The observed $y$ value has a positive prior predictive density,
    i.e., $\E{\psi(y|X,Z,\Theta)}>0$.
  \end{enumerate}
\end{assum}

In what follows, we will denote by $\cmspace^d$
the space of absolutely continuous $x\colon\tspace\to\reals^d$ with
square integrable weak derivatives $\dot x$.
For $x\in\reals^d$, $\norm{x}$ will denote its Euclidean norm.
Furthermore, $\snorm{\cdot}$ will denote the supremum norm of
continuous functions from $\tspace$ to $\reals^d$:
\begin{align}
  \label{eq:supremum_norm_defn}
  \snorm{w} &:= \textstyle\sup_{t\in\tspace} \norm{w(t)}.
\end{align}
The divergence of a vector field function $f$, with respect to a variable $x$
is denoted $\diverg_x f$, i.e.,
$\diverg_x f = \sum_k \frac{\partial f_k}{\partial x_k}$.

The theorem below is the application of Defn.~\ref{th:fict_dens}
to the state paths and parameters of systems defined by SDEs.
Its detailed proof is presented in a previous work
\citep[Sec.~2.2]{dutra2014thesis}.
The proof uses a generalized version stochastic Stokes' theorem and
follows roughly the same steps as the derivation of the Onsager--Machlup 
functional by \citet[Thm.~2]{capitaine2000omf}.
The extension of the Onsager--Machlup functional to the unknown parameter
case is nontrivial and an original result by the authors.
The assumptions on the model parametrization under which the theorem is valid 
also give insight into the applicability of the resulting estimator.

\begin{thm}[{\citealt[Thm.~2.22]{dutra2014thesis}}]
  \label{th:main_thm}
  When Assum.~\ref{th:assum} holds, then the joint fictitious density
  $\rho$ of $X$ under the $\snorm{\cdot}$ norm, and 
  $Z_0$ and $\Theta$ under the Euclidean norm, conditioned on $Y=y$,
  is given by
  \begin{align}
    \label{eq:post_joint_fict_dens}
    \rho(x,z_0,\theta|y) =
    \frac{
      \psi(y|x,z,\theta)\pi(x(0), z_0, \theta)\exp\big(J(x,z_0, \theta)\big)
    }{
      \E{\psi(y|X,Z,\Theta)}
    },
  \end{align}
  for all $x\in\cmspace^\nx$, $z_0\in\reals^\nz$ and $\theta\in\reals^\npar$,
  where $J$ is the Onsager--Machlup functional, defined by
  \begin{multline}
    \label{eq:om_func_defn}
    \textstyle
    J(x,z_0, \theta) :=
    -\frac12 \int_0^T \diverg_x f(t, x_t, z_t, \theta)\,\dd t\\
    \textstyle
    -\frac12 \int_0^T \normnm{
      G\inv[\dot x_t - f(t, x_t, z_t, \theta)]
    }^2\dd t,
  \end{multline}
  and $z\in\cmspace^\nz$ is the unique solution to the initial
  value problem (IVP)
  \begin{align}
    \label{eq:zeta_ivp}
    \dot z(t) & = h\fparen{t, x(t), z(t), \theta}, &
    z(0)& = z_0.
  \end{align}
\end{thm}

MAP state-path estimation is often done using the Euler-discretized
SDE or by omitting the drift divergence term in \eqref{eq:om_func_defn}
\citep{aravkin2011llr, aravkin2012rtf, varziri2008pse, karimi2014mlm}.
The theorem below, which is also proved by \citet[Sec.~2.3]{dutra2014thesis},
shows that by doing so the fictitious density 
\emph{of the associated noise path} is obtained instead.

\begin{thm}[{\citealt[Thm.~2.26]{dutra2014thesis}}]
  \label{th:noise_map_thm}
  When Assum.~\ref{th:assum} holds, then the joint fictitious density
  $\rho_{\mathrm e}$ of \ $W$ under the $\snorm{\cdot}$ norm, and $X_0$,
  $Z_0$ and $\Theta$ under the Euclidean norm, conditioned on $Y=y$,
  is given by%
  \footnote{
    \sloppy
    We note that $\E{\psi(y|X,Z,\Theta)}$ is missing from Thm.~2.26 of
    \citet{dutra2014thesis} due to a typographical error.
  }
  \begin{align}
    \label{eq:noise_post_joint_fict_dens}
    \rho_{\mathrm e}(w,x_0,z_0,\theta|y) =
    \frac{
      \psi(y|x,z,\theta)\pi(x_0, z_0, \theta)
      \exp\big(J_{\mathrm e}(x,z,\theta)\big)
    }{
      \E{\psi(y|X,Z,\Theta)}
    },
  \end{align}
  for all $x_0\in\reals^\nx$, $z_0\in\reals^\nz$, $\theta\in\reals^\npar$
  and $w\in\cmspace^\nx$ such that $w(0)=0$, where $J_{\mathrm e}$ is 
  the energy functional, defined by
  \begin{equation}
    \label{eq:energy_func_defn}
    \textstyle
    J_{\mathrm e}(x,z, \theta) :=
    -\frac12 \int_0^T \normnm{
      G\inv[\dot x_t - f(t, x_t, z_t, \theta)]
    }^2\dd t,
  \end{equation}
  and $x\in\cmspace^\nx$ and $z\in\cmspace^\nz$ are the
  unique solutions to the IVPs
  \mathtoolsset{showonlyrefs=false}%
  \begin{subequations}
    \label{eq:xz_ivp}
    \begin{align}
      \label{eq:x_ivp}
      \dot x(t) & = f\fparen{t, x(t), z(t), \theta} + G\dot w(t), &
      x(0)& = x_0, \\
      \label{eq:z_ivp}
      \dot z(t) & = h\fparen{t, x(t), z(t), \theta}, &
      z(0)& = z_0.
    \end{align}
  \end{subequations}
  \mathtoolsset{showonlyrefs=true}%
\end{thm}

We note that the divergence of a vector field quantifies the expansion or 
contraction of an infinitesimal volume around a given point.
This means that the probability of the neighbourhood of a state-path and
parameter vector depends not only on the probability of the neighbourhood of the
associated noise path, initial states and parameters,
it also depends on the amplification or attenuation by the drift $f$ of the
pertubations in $X$ due to the process noise around the chosen path.

From the definition of MAP \citep[Defns.~1-2]{dutra2014map} and of
fictitious density, it follows from Thm.~\ref{th:main_thm} that the joint
MAP state path and parameter estimates are the maxima of the fictitious density
\eqref{eq:post_joint_fict_dens}.
For the purpose of maximizing the fictitious densities, the prior predictive
density---the denominator of \eqref{eq:post_joint_fict_dens}---can be omitted
as it is constant for each experiment and its observed $y$.
Furthermore, the logarithm of the fictitious density can be used as the merit
function as it leads to a better-conditioned optimization problem.
Hence the joint MAP estimate of $X$, $Z_0$ and $\Theta$ is the solution to the
optimal control problem below:
\begin{equation}
  \label{eq:jme_problem}
  \begin{aligned}
    \operatorname*{maximize}_{
      x\in\cmspace^\nx\;\,z\in\cmspace^\nz\;\,\theta\in\reals^\npar
    } &&\;&
    \ell(x,z,\theta,y)\\
    \operatorname{subject\,to}\;\;\quad&&&
    \dot z(t) = h\fparen{t, x(t), z(t), \theta},
  \end{aligned}
\end{equation}
where the merit function $\ell$ is given by
\begin{multline}
  \label{eq:logpost}
  \ell(x,z,\theta,y) = 
  \ln\psi(y|x,z,\theta)\\
  +\ln\pi(x_0, z_0, \theta)
  \textstyle
  -\frac12 \int_0^T \diverg_x f(t, x_t, z_t, \theta)\,\dd t\\
  \textstyle
  -\frac12 \int_0^T \normnm{
    G\inv[\dot x_t - f(t, x_t, z_t, \theta)]
  }^2\,\dd t.
\end{multline}

Thm.~\ref{th:noise_map_thm} implies that if the drift divergence is omitted
from \eqref{eq:logpost}, then the state path
associated with the joint MAP noise path, initial states and parameters are
obtained.
Following \citet{dutra2014map}, we denote this the
joint minimum-energy state path and parameter estimator (MEE),
which is solution to the optimal control problem below:
\begin{equation}
  \label{eq:mee_problem}
  \begin{aligned}
    \operatorname*{maximize}_{
      x\in\cmspace^\nx\;\,z\in\cmspace^\nz\;\,\theta\in\reals^\npar
    } &&\;&
    \lenerg(x,z,\theta,y)\\
    \operatorname{subject\,to}\;\;\quad&&&
    \dot z(t) = h\fparen{t, x(t), z(t), \theta},
  \end{aligned}
\end{equation}
where the merit function $\ell_e$ is given by
\begin{multline}
  \label{eq:energy_logpost}
  \lenerg(x,z,\theta,y) = 
  \ln\psi(y|x,z,\theta)
  +\ln\pi(x_0, z_0, \theta)\\
  \textstyle
  -\frac12 \int_0^T \normnm{
    G\inv[\dot x_t - f(t, x_t, z_t, \theta)]
  }^2\,\dd t.
\end{multline}
The estimation problems \eqref{eq:jme_problem} and \eqref{eq:mee_problem}
can then be solved with standard optimal control techniques such as direct 
collocation \citep[cf.][and the references therein]{betts2010pmo}.

\section{Simulated examples}
\label{sec:sim}
To illustrate the proposed estimator, we apply it on simulated data of the
Duffing oscillator, a benchmark model for modeling nonlinear dynamics and
chaos \citep{aguirre2009mnd} and state estimation in SDEs
\citep{ghosh2008sis, khalil2009nfc, namdeo2007nsd}.
We note that the standard Duffing oscillator does not satisfy
Assum.~\ref{th:assum}.
However, the equations are of no interest too far from the origin, so an
equivalent system satisfying the assumptions can be obtained by multiplying
the drift function by a $C^2$ bounding function which is equal to unity
in the region of interest and is zero for large values of the states.

The system dynamics is given by
\begin{align}
  \dd X_t &= 
  \phi(Z_t)[-AZ^3_t -BZ_t - DX_t + \gamma\cos t]\dd t + \sigmad\dd W_t\\
  \dd Z_t &= X_t \, \dd t,
\end{align}
where $A$, $B$ and $D$ are parameters considered unknown, to be estimated, and 
$\gamma$ and $\sigmad$ are parameters considered known, whose nominal values
are used in the estimation, and $\phi\colon\reals\to\reals$ is the
smootherstep bounding function:
\begin{align}
  \phi(z) &=
  \begin{cases}
    1 & |z|\leq 1000\\
    0 & |z|\geq 1001\\
    \eta(1001 - |z|) & \text{otherwise},
  \end{cases}\\
  \eta(\varepsilon) &= 6 \varepsilon^5 - 15\varepsilon^4 + 10 \varepsilon^3.
\end{align}
We note that in no case the simulated $|Z|$ or estimated $|z|$
were larger than \np{1000}, so the results would be unchanged had the
bounding function not been used.

Noninformative Gaussian priors were chosen for the drift-function parameters,
with a large standard deviation:
\begin{align}
  A&\sim\normald{0, \sigmap^2}, &
  B&\sim\normald{0, \sigmap^2}, &
  D&\sim\normald{0, \sigmap^2}.
\end{align}

The system was simulated using the strong explicit order 1.5 scheme 
\citep[Sec.~11.2]{kloeden1992nss} with a time step of 0.005.
The initial states were sampled from the independent normal distributions
$X_0\sim\normald{0,\sigmai^2}$ and $Z_0\sim\normald{0,\sigmai^2}$ and
all realizations of the simulations were performed with the parameters
at their nominal values, shown in Table~\ref{tab:nominal}.
A total of 100 Monte Carlo simulations were performed for each experiment.

\begin{table}
  \centering
  \caption{Nominal parameter values of the Duffing oscillator experiments.}
  \label{tab:nominal}
  \begin{tabular}{cccccccccc}
    \hline
    $A$ & $B$ & $D$ & $\gamma$ & $\sigmad$ & $\sigmai$ & $\sigmap$ & $r$ & $s$ &
    $\ts$ \\
    1.0 & -1.0 & 0.2 & 0.3 & 0.1 & 0.4 & 10 & 1.1 & 10 & 
    0.1 \\
    \hline
  \end{tabular}
\end{table}

The JME \eqref{eq:jme_problem} and MEE \eqref{eq:mee_problem} were applied to
the simulated data.
These estimators were implemented using the same optimal control techniques
used in \citep[Sec.~9]{dutra2014map}.
The estimation was transcribed to a nonlinear programming (NLP) problem
using a third-order Legendre--Gauss--Lobatto direct collocation method
equivalent to the Hermite--Simpson method \citep[Sec.~4.5]{betts2010pmo}.
The resulting NLP was then solved using the IPOPT solver of 
\citet{wachter2006iip}.

The results were then compared with the MAP parameter estimates obtained
with the PEM \citep[cf.][]{kristensen2004pes}, which was 
implemented using the unscented Kalman filter (UKF) with the 
SDE prediction step of \citet{arasaratnam2010ckf}.
The unscented Kalman smoother (UKS) with the backward-correction step of
\citet{sarkka2008urt} was then used to obtain the optimal state-path
associated with the PEM-estimated MAP parameters.
To be more favorable with the PEM and avoid local minima, the nominal
parameter values were used as its optimization's starting point.

Noise-corrupted discrete-time measurements of the $Z$ state were taken
with period $\ts$.
Two experiments with different measurement noise distributions were done
to investigate different aspects of the proposed estimator.
Each is detailed in the following subsections.

\subsection{Measurements with Gaussian noise}
\label{sec:gaussian_noise}

In the first example, Gaussian measurement noise was used.
In this way, the UKF, UKS and PEM are applicable and can serve as benchmarks.
Each measured value $Y_k$ was drawn independently from
\begin{equation}
  \label{eq:gauss_meas}
  Y_k |X,Z,\Theta \sim \normald{Z_{k\ts}, \Sigmay^2},
  \qquad k=0,\dotsc,N,
\end{equation}
where the standard deviation $\Sigmay$ is a parameter considered unknown,
to be estimated.
The nominal value $\Sigmay=0.1$ was used to generate the data.
For the estimation of $\Sigmay$, the gamma distribution with shape $r$ and
scale $s$ was chosen as its prior.
The total experiment length $T=200$ was used.

The initial guess (optimization starting point) for the $z$ path 
of the MEE and JME was obtained by doing
a least-squares spline fit to the measured data.
The initial guess for the $x$ path was then the derivative of this spline.
For the initial guess of the drift-function parameters a least-squares
regression was performed using the second derivative of the spline and
the guess of the $x$ and $z$ paths.
Finally, the initial guess of the $\sigmay$ parameter was the sample standard
deviation of $Y_{k\ts} - z_{k\ts}$, using the initial $z$ guess.
\begin{figure}
  \centering
  \inputtikzpicture{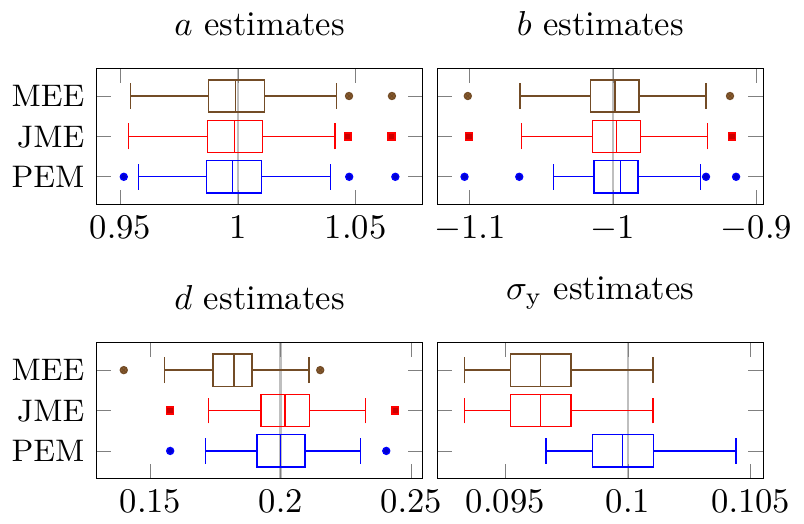}
  \caption{%
    Boxplots of estimated parameters of the Duffing oscillator with Gaussian
    measurement noise.
    The nominal values are marked with gray vertical lines.
  }
  \label{fig:duffing_popt}
\end{figure}

\begin{figure}
  \centering
  \inputtikzpicture{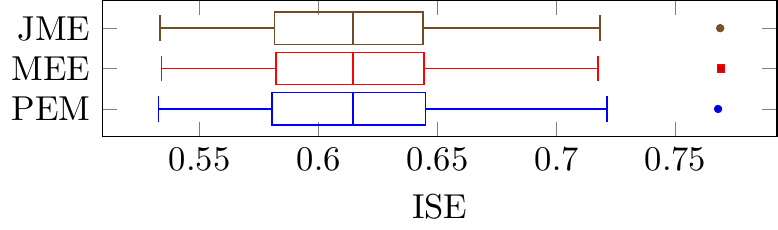}
  \caption{%
    Boxplot of the state-path estimation error of the Duffing oscillator
    with Gaussian measurement noise.
  }
  \label{fig:duffing_ise}
\end{figure}

In Fig.~\ref{fig:duffing_popt}, boxplots of the parameter estimates of
the three estimators are shown.
It can be seen that the distribution of the estimates of the drift-function
parameters by the JME and PEM is very similar.
The MEE, on the other hand, is biased for the $d$ parameter, which
represents the damping of the system.
This is because the drift divergence $\diverg_x f(t, x, z, \theta) = -d$.
Consequently, higher values for the $d$ parameter imply that the system
strongly attenuates the fluctuations due to the process noise and increase
the probability of small balls around the state path.
The MEE, which does not take this into account, underestimates this parameter.
It should also be noted that both the MEE and JME are biased for the $\sigmay$
parameter, as lower values of this parameter also increase the
probability of small balls around state paths.

To evaluate the state-path estimation error, the integrated square error (ISE)
metric was used, given by
\begin{equation}
  \textstyle
  \operatorname{ISE} = \int_0^T [(X_t - x_t)^2 + (Z_t - z_t)^2]\,\dd t,
\end{equation}
where $X$ and $Z$ are the simulated paths and $x$ and $z$ their estimates.
Boxplots of the ISE distribution of the three estimators is shown in
Fig.~\ref{fig:duffing_ise}, in which we see that their error is comparable.

\subsection{Measurements with outlier noise}

\begin{figure}
  \centering
  \inputtikzpicture{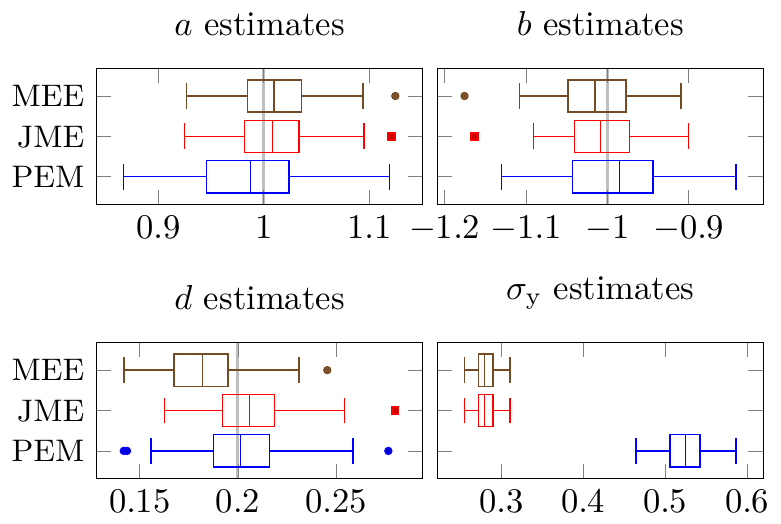}
  \caption{%
    Boxplots of estimated parameters of the Duffing oscillator with outlier
    measurement noise.
    The nominal values are marked with gray vertical lines.
  }
  \label{fig:dr_popt}
\end{figure}

\begin{figure}
  \centering
  \inputtikzpicture{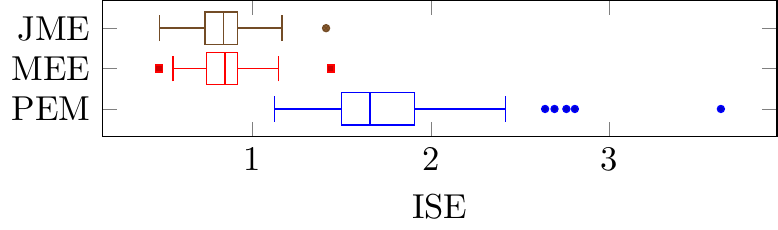}
  \caption{%
    Boxplot of the state-path estimation error of the Duffing oscillator 
    with outlier measurement noise.
  }
  \label{fig:dr_ise}
\end{figure}

For the second example, measurement noise with outliers was used, to
illustrate the advantages of the proposed estimator for robust system
identification and smoothing.
Each measured value $Y_k$ was drawn independently from the Gaussian
mixture distribution
\begin{equation}
  \label{eq:outl_meas}
  Y_k |X,Z,\Theta \sim \po \normald{Z_{k\ts}, \sigmao^2}
  + (1 - \po) \normald{Z_{k\ts}, \sigmar^2},
\end{equation}
where $\po=0.25$ is the outlier probability, $\sigmao=1$ is the outliers'
standard deviation, and $\sigmar=0.2$ is the regular measurements' standard
deviation.
The total experiment length $T=100$ was used.
This experiment is similar to those of \citet{aravkin2012rtf} and
\citet{dutra2014map}.

For estimation with the MEE and JME, a measurement model different from the
one used to generate the data was used in the estimation.
Student's $t$-distribution with 4 degrees of freedom and unknown scale $\Sigmay$
was used as the measurement likelihood, with the following expression
for its log-likelihood:
\begin{equation}
  \label{eq:t_log_likelihood}
  \textstyle
  \ln\psi(y|x,z,\theta) = -\sum_{k=0}^N
  \left[
  \frac52\ln\big(1 + \frac{(y_k - z_{k\ts})^2}{4\sigmay^2}\big) + \ln\sigmay
  \right],
\end{equation}
where the constant terms have been omitted as they do not influence the
location of maxima.
For the estimation of $\Sigmay$, the gamma distribution with shape $r$ and
scale $s$ was chosen as its prior, as in the previous experiment.
The optimization starting point for the MEE and JME was also the same as in
Sec.~\ref{sec:gaussian_noise}.

Boxplots of the parameter estimates by the three estimators are shown
in Fig.~\ref{fig:dr_popt}.
It can be seen that, as in the previous example, the MEE is biased for
the $d$ parameter.
In this experiment, however, the PEM also presents a larger interquartile range
than the JME.
The difference between the PEM and JME is clearer on
the state-path estimation error, shown in Fig.~\ref{fig:dr_ise}.
The JME and MEE achieve lower errors because the $t$-distribution fits
the heavy-tailed distribution used to generate the measurements better than
the Gaussian distribution which underlies the PEM and UKS.
We note that no nominal value is marked for $\sigmay$ as a measurement model
different from the one used to generate the data was used for the estimation.
Furthermore, this parameter has a different meaning for the PEM and for the 
JME--MEE.

\section{Conclusions and future work}
\label{sec:conclusions}
In this paper, we extended the estimator of \citet{dutra2014map} for joint
MAP and minimum-energy state-path and parameter estimation.
The minimum-energy estimator is proved to estimate the joint MAP noise paths,
initial states, and parameters.
The difference between the MAP and minimum-energy estimators is that the
former takes into account the amplification or attenuation by the system of
the pertubations due to the process noise.
It does so by including the drift-divergence integral in its log-posterior.
Simulated experiments show that this is important for preventing the biasing
of parameters related to the system's energy-loss rate.

These results suggest that the Onsager--Machlup functional should be used
instead of the energy functional in approximate MAP or maximum likelihood
estimation by coestimating the state path as nuisance parameter
\citep[cf.][]{varziri2008pse}.
It also indicates that the Onsager--Machlup functional should be used
in MAP parameter estimators which marginalize the state-path using the
Laplace approximation \citep{karimi2014mlm}, which we intend to investigate
in future work.

The JME and MEE were also compared with the PEM, which occupies a similar
niche of applications as its parameters are estimated as single vector and not
as a time-varying signal.
In the example with Gaussian noise, the drift-function parameter estimates
by the JME were comparable to those by the PEM.
However, one advantage of the JME is that it is applicable to systems with
a wider range of measurement distributions.
In the experiment with non-Gaussian measurements, the JME obtained lower errors
than the PEM.
One serious disadvantage of the JME, however, is that it cannot be used to
estimate the diffusion matrix, unlike the PEM.

\bibliographystyle{plainnat}
\bibliography{bibtex-compressed}

\end{document}